# A BAYESIAN HIERARCHICAL MODELING APPROACH TO DIETARY ASSESSMENT VIA FOOD FREQUENCY

## Abstract

Previous likelihood-based linear modeling of nutritional data has been limited by the availability of software that allows flexible error structures in the data. We demonstrate the use of a Bayesian modeling approach to the analysis of such data. Our goal is to model the relationship between the energy intake derived from Food Frequency Questionnaires (FFQs) and the energy expenditure estimated from the doubly labeled water method. We consider models with different distributions for the FFQ energy intake. The models include previously identified covariates describing social desirability, education and their possible interaction that are felt to impact the reported FFQ. The models also include random effects to account for subject specific random variation (frailty) and also to account for the complex patterns of measurement error inherent in these data. Issues arising within the work relate both to the selection of relevant linear and non-linear models, the use of random effects, and the relevance of goodness-of-fit criteria such as DIC and PPL in assessing the most appropriate model.

## 1 Introduction

In dietary nutrition studies it is of great importance to measure accurately food intake and to assess the relation of food intake to energy expenditure and intake. To this end the food frequency questionnaire (FFQ) is a standard tool for dietary assessment of 'normal' intake (over a fixed period) and is often administered at the start and sometimes the end of a nutritional study. FFQs consist of a list of foods and also ask about the usual portion size the individual consumes. This information is converted to energy intake (EI). The FFQ is self-reported and is known to be associated with both systematic and random errors ((Rosner and Gore2001;Carroll et al.1998;Rosner and Gore2001).

In the analysis of such data there is frequently a need to relate food intake to other measured covariables describing individuals' background and lifestyle. Such variables can be, for example, education level (or years in education) which may relate to the nutritional habits and physical activity, or psycho-social constructs such as social desirability, which relates to individual's self image. Accurate reporting of food intake can relate closely to these variables. Hence it is useful to evaluate these covariates in relation to intake to assess whether there are variations explained by such variables.

Using doubly labeled water, one can obtain a measure of Total Energy Expenditure (TEE) free of self-report bias. Under certain conditions of energy balance, the energy expenditure can be assumed to be equal to the energy intake. Thus, the relation between the TEE measured by the doubly labeled water and the EI estimated from the FFQ can provide insight into the measurement error and biases associated with the self-reporting.

In the following we examine a range of novel Bayesian hierarchical models for the relation between total EI measured by FFQ (Kcal/day) and the objective assessment of intake based on doubly labeled water (DLW: Kcal/day). We develop models which admit a variety of types of error in the measured DLW and in the relation of FFQ to DLW. In the following section we outline a dataset which has been derived from a study of energy balance in women (the ENERGY study). Issues associated with the analysis of these data prompted our interest in modeling the FFQ relationship, and we will use the data to help illustrate our methodology. The purpose of our study is twofold. Primarily, we want to assess the model adequacy for nutritional data using conventional assumptions (such as normal error structures for overall error and random effects). Secondly, we were interested in assessing the relation between a measure of energy based on FFQ and a measure based on DLW in the presence of specific covariates . The measures included are described in Section 2.1. The rest of the paper is organized as follows. In the next section the dataset and measures available are described. In section 3 basic models for the relation between FFQ and DLW are discussed, and these involve linear mixed models which include measurement error. In section 4 we described further the hierarchical components of the models including prior distributions for random effect variances, and hyperprior distributions. Section 5, 6 and 7 deal with model fitting, evaluation metrics and preliminary results for the linear mixed model. Extensions to generalized linear mixed models both in terms of alternative error structures and random effect prior distributions are discussed in section 8.

Finally the results of the model selection are discussed in section 9, both in the context of evaluation of the linear mixed model, and in terms of the effect of important covariates.

## 2 Data Set

A complete description of study methods and data collection can be found elsewhere (Hebert et al.2002). Briefly, 81 women in Worcester, MA aged 40 to 65 years were recruited into the study between June and October 1997. Eligible women were free of major medical conditions, not on steroid therapy (including inhaled steroids), and had a body mass under 91 kg (200 lbs). Women could be taking estrogen replacement or birth control hormones. They were excluded if they were on a special diet to lose or gain weight, if they would not consent to maintaining their current physical activity and dietary habits during the study, or if they would not be available to be reached by telephone on any day during the study period. Premenopausal subjects were scheduled to enter the study the week after they completed menstruation. All procedures were approved by the Institutional Review Board of the University of Massachusetts Medical School.

## 2.1 Measures

On each participant a range of measures were recorded. Herein, we examine the following variables:

<u>Doubly labeled water (DLW)</u>: Following oral administration of DLW ($^2H_2^{18}O$), the $^{18}O$ label is eliminated from the body as both carbon dioxide and water, and the deuterium label is excreted exclusively as water. Thus, the difference between the urinary elimination rates of $^{18}O$ and deuterium provides a measure of carbon dioxide production ($rCO_2$). Total Energy Expenditure (TEE) can be calculated from $rCO_2$ using an estimate of respiratory quotient (RQ) (7, 8).

<u>Dietary Energy Intake (EI)</u>: The FFQ employed in this study is identical to that used in the Women's Health Initiative (WHI) (6). Based on the instrument developed by Block and colleagues at NCI (9), the FFQ is a 145-item questionnaire consisting of 19 food behavior questions, 114 food and 8 beverage frequency questions, and 4 summary questions. The instrument was designed to assess >90% of energy and fat intake in the U.S. population as a whole, though this proportion would vary between population subgroups. The database from which nutrients were computed for the FFQ also is derived from the NDS. The FFQ was administered twice, once just prior to collection of the first urine sample (day 0) and the other just prior to the end of participation (day 13). Because the vast majority of data in most epidemiologic studies are derived from the first-use FFQ (even in studies where >1 FFQ is used, subsequent ones are often administered years apart), we used estimates based on the day 0 administration here.

Psychosocial Measures:  *Social desirability* was measured using the Marlowe-Crowne Social Desirability Scale (MCSD) (11) whose 33 true / false questions quantify a tendency to avoid criticism and to defend one's social image in a testing situation, portraying oneself as conforming to societal expectations (e.g., "I never hesitate to go out of the way to help someone in trouble").  This scale has been used in the validation and editing of hundreds of psychological questionnaires over the past 35 years (12-17). It has a test-retest reliability coefficient of .89, showing relatively little change over time, and an internal consistency coefficient of .88 (11).

Education (Edu): The *education level* was recorded as having a college degree or less than college degree.

## 3  Models for FFQ and DLW dependence

In what follows, we assume that the outcome variable is the total FFQ energy value per individual. This we denote for *n* responses as $\boldsymbol{y} = \{y_1,....,y_n\}$. The independent variable that we want to relate to FFQ outcome is doubly labeled water (DLW). We denote this by $\boldsymbol{x}_1 : \{x_{11},....,x_{in}\}$. We are also particularly interested in the relation between FFQ and DLW in the presence of specific covariates.  The covariates of interest are social desirability and education level as described in section 2. Social desirability is measured on a continuous scale, whereas education level is a binary factor with threshold at college degree. These two additional covariates are denoted $\boldsymbol{x}_2, \boldsymbol{x}_3$. Consider first a generic generalized linear model for the relationship of the form:

$$\boldsymbol{y} \sim f(\boldsymbol{\mu})$$

where the expectation of $\boldsymbol{y}$ is assumed to be $E(\boldsymbol{y}) = \boldsymbol{\mu}$. Our basic starting point is a form that is conventionally assumed within dietary studies ((Carroll et al.1998;Kipnis et al.2001), (Kipnis et al.2003), (Berry, Carroll, and Ruppert2002)), that of the Gaussian error model. Hence we assume that $\mathbf{y} = \boldsymbol{\mu} + \mathbf{e}$ where the mean parameter $\boldsymbol{\mu}$ is parameterized with ( a linear) function of covariates. Often this linear model is assumed once a log transformation is made of both dependent and independent variables (see, e.g. (Carroll et al.2006)). We will return to this in a later section.

Define first, for each observation, a Gaussian mixed model:

$$(y_i \mid \mu_i, e_i) \sim N(\mu_i, e_i)$$
$$g(\mu_i) = x^T_{\ i}\beta + z^T_{\ i}\xi. \qquad (1)$$

Here we define a link function $g(\mu_i)$ which relates the mean level of the outcome to a linear predictor $x^T_{\ i}\beta + z^T_{\ i}\xi$. The identity link function is assumed. Note that $x^T_{\ i}$ is the vector of fixed covariate observations for the $i$ th observation from the design matrix defined for our example by $\mathbf{x} = \begin{Bmatrix} 1 & x_{11} & x_{21} & x_{31} \\ . & . & . & . \\ . & . & . & . \\ 1 & x_{1n} & x_{2n} & x_{3n} \end{Bmatrix}$, and $\beta$ is the associated parameter vector, while $z^T_{\ i}$ is the random effect vector for the $i$ th observation whose length corresponds to the number of individual-level effects defined in unit vector $\xi$. In matrix notation this leads to the model

$$\mathbf{y} \sim N(\mathbf{x}\beta + \mathbf{z}\xi, \tau I_n)$$

where $\Sigma_{\mathbf{\epsilon}} = \tau I_n$ is the covariance matrix of the overall model error ($\mathbf{e}$).

In the first model, FFQ is an outcome, and we consider DLW as a fixed covariate with other fixed covariates being the social desirability ($x_2$) and the education level ($x_3$). The standard link function assumed for Gaussian error is the identity link, hence the model $\mu_i = x^T_{\ i}\beta + z^T_{\ i}\xi$ is commonly assumed. We start with a Gaussian model with no random effect component and g(.) as the identity function, i.e.

$$\mu_i = x^T_{\ i}\beta \qquad (1a).$$

 This model assumes that the covariates are fixed and that there is no measurement or other error in the measured covariates. As we assume a Bayesian modeling paradigm we do allow the regression parameters to have prior distributions. Prior distributions for parameters and effects are discussed in a later section (Section 4.1).

Model (1a) above is a standard Bayesian linear model, with Gaussian likelihood. It is commonly assumed however that various types of error may contaminate the measurement of covariates, particularly in nutritional food intake studies. This is true for DLW as for other measures, and so we also consider an extension to models with additional error, of various forms. Here we particularly focus on error in DLW as the relation of that to FFQ is the primary focus of the analysis.

 The first of these models extends the Gaussian model by adding random effect components. In the generic model (1), we have included a vector of random effects as additive components. This is a simplification of a more general specification. In our model extension we admit both additive effects in the linear predictor and also additive effects in the measured covariates. For example with a single random effect $z_{i\ i}$ where $g(.)$ is still the identity function, the linear predictor becomes

$$\mu_i = x^T_{\ i}\beta + z_i . \qquad (1b)$$

Assumption of a Gaussian error model leads to

$$y_i = x^T_{\ i}\beta + z_i + e_i. \qquad (1c)$$

This model is a simple variance component model with two components based on ($z_i, e_i$). As such, this model allows for extra variability related to the observation unit level but does not specifically relate this variability to error in measurement in covariates. Of course ($z_i, e_i$) are not identified in the likelihood without further prior distributional assumptions. Hence they can be Bayesianly identified (Gelfand and Sahu 1999). A more direct parameterization of measurement error would be the following:

$$\mu_i = \beta_0 + \beta_1 x_i^*$$
$$where\ x_i^* = x_i + \varepsilon_i .$$

Here $x_i$ is the observed covariate which is assumed to relate to the true covariate ($x_i^*$) by addition of an error component ($\varepsilon_i$). Note that in this formulation the error contributes also to the relationship with the outcome as $\beta_0 + \beta_1 x_i^* = \beta_0 + \beta_1 x_i + \beta_1 \varepsilon_i$. As this more closely identifies the error related to a particular covariate we adopt this as a first stage measurement error model. When applied to measurement error in DLW, this model yields the following full specification

$$\mu_i = \beta_0 + \beta_1 \left( x_{1i} + \varepsilon_i \right) + \beta_2 x_{2i} + \beta_3 x_{3i} + \beta_4 \left( x_{2i} * x_{3i} \right) \qquad (2).$$

Here it is assumed that any error in $x_{2i}, x_{3i}$ is negligible compared to the error in DLW.

# 4 Further Model Components

Hierarchical Bayesian models consist of levels of variation described by (prior) distributions. In our measurement error and random effect models we must specify the prior distributions for all random components.

## 4.1 Prior distributions

In general, conventional specifications are assumed for prior distributions. For regression parameters we make the usual zero-mean Gaussian assumption, i.e. $\beta_* : N(0, \tau_{\beta_*})$ where the $*$ denotes a generic regression parameter. This specification is reasonably non-informative, particularly if $\tau_\beta$ is taken to be very large. A more controversial specification for the variance parameter ($\tau_{\beta_*}$) is the use of gamma or inverse gamma hyper-prior distributions. While gamma prior distributions are commonly used for precision parameters, this has recently been criticized by Gelman (Gelman2006) who shows that such prior distributions can be highly informative. Instead he recommends the use of uniform distributions on large positive range for standard deviation hyper-priors. Uniform distributions as priors for variances can lead to improper posterior distributions, however, and so must be applied carefully.

To fully construct model (2) we assigned the following hierarchical structure:

$$[y_i \mid \mu_i, \tau_y] \sim N(\mu_i, \tau_y)$$
$$\mu_i = \beta_0 + \beta_1 \left(x_{1i} + \varepsilon_i\right) + \beta_2 x_{2i} + \beta_3 x_{3i} + \beta_4 \left(x_{2i} * x_{3i}\right).$$

At the next level we have:

$$[\beta_* \mid \tau_{\beta_*}] \sim N(0, \tau_{\beta_*}),$$

and for the random effects in model (2) we assume: $[\varepsilon_i \mid \tau_\varepsilon] \sim N(0, \tau_\varepsilon)$ . The next level of the hierarchy requires the specification of hyper-prior distributions for the variance parameters of the regression parameters and the random effects (consisting of measurement error ($\varepsilon$) and model error ($e$)). The distributions assumed are:

$$s_* \sim U(0, B_*) \text{ where } s_* = \sqrt{\tau_{\beta_*}},$$

following Gelman's suggestion, and

$$s_\varepsilon \sim U(0, B_\varepsilon) \text{ where } s_\varepsilon = \sqrt{\tau_\varepsilon},$$

and

$$s_y \sim U(0, B_y) \text{ where } s_y = \sqrt{\tau_y}.$$

The B parameters are fixed a priori. Here we assume that $B_* = 1000$, $B_e = 1000$, and $B_\varepsilon = 250$. Thus the variance of the random effect associated with the measurement error in the DLW is smaller than the model error.

With the above complete hierarchy we can consider the posterior distribution and associated sampling.

## 5 Posterior Distributions and Model Fitting

Model (2) is a Gaussian likelihood model with Gaussian prior distributions and uniform hyper-priors. The posterior distribution is given by

$$P_0(\beta_*, \boldsymbol{\varepsilon}, s^2_{\varepsilon}, s^2_{\beta_*}, s^2_{\varepsilon} \mid \mathbf{y}, \mathbf{x}) \propto N(\mathbf{y} \mid \mathbf{x}, \beta_*, \boldsymbol{\varepsilon}, s^2_{\varepsilon})[\beta_* \mid s^2_{\beta_*}][\boldsymbol{\varepsilon} \mid s^2_{\varepsilon}][s_{\beta_*} \mid B_*][s_{\varepsilon} \mid B_{\varepsilon}][s_{\varepsilon} \mid B_{\varepsilon}]$$

where $N(\mathbf{y} \mid \mathbf{x}, \beta_*, \boldsymbol{\varepsilon}, s^2_{\varepsilon})$ denotes a normal likelihood for data $\mathbf{y}$ given $\mathbf{x}, \beta_*, \boldsymbol{\varepsilon}, s^2_{\varepsilon}$.

This form of the model has a log concave posterior distribution and so adaptive rejection sampling can be used to provide samples form the conditional distribution of all parameters. Hence, it is straightforward to use Gibbs sampling for this model. In fact, WinBUGS has been used to provide McMC samples from the posterior distribution (Gamerman2002). As with all McMC sampling, convergence must be checked before taking final sample values. We used 100,000 iterations for burn-in for each model fit. We ran multiple chains to assess the convergence of the process, based on the Brooks-Gelman-Rubin (BGR) statistic (Brooks and Gelman1998). As noted in the next section we examined a range of goodness-of-fit criteria for models. One such measure is the posterior predictive loss (PPL). For convergence checking we examined both streams of single parameters and also the PPL using BGR.  To reduce autocorrelation we thinned the final sample. To achieve this we took every 50th iteration up to a total sample of 2000 observations after burn-in.

## 6 Model Selection

In our analysis of the Hebert et al. dataset we examined a variety of error structures and models. Model (2) is used as a starting point for our model selection process.  We used two main criteria to assess model adequacy. First, the deviance information criterion (DIC) proposed by Spiegelhalter (Spiegelhalter et al.2002) is used to compare different models. This is defined as $DIC = \bar{D} + pD = \bar{D}(\boldsymbol{\theta}) + (\bar{D}(\boldsymbol{\theta}) - D(\bar{\boldsymbol{\theta}}))$ where $\bar{D}(\boldsymbol{\theta})$ and $D(\bar{\boldsymbol{\theta}})$ are the average deviance of parameters and deviance of the posterior average parameters respectively.

The model with the smallest DIC indicates the best model. If the posterior distribution is very skewed or even multimodal, the DIC may not be appropriate. Use of transformations can dramatically affect DIC (in fact log transformations of the dependent variable can produce negative estimates of effective number of parameters ($pD$)).

To aid in model selection, we examined another criterion: the posterior predictive loss as measured by the mean square prediction error (MSPE) under the squared error loss function as described in Gelfand and Ghosh (1998). The MSPE is the mean squared difference between the observed and the predicted values of the outcome variable. Thus, the model which results in predicted values closest to the observed values will produce the lowest MSPE. This appears to be more robust to scale change and to transformation.

## 7 Initial Results

Model (2) above was fitted to the Hebert et al. data in three different forms: linear regression model with no random effects as in (1a); linear regression model with random effects as in (1c); and linear regression model with random effects inside the regression on DLW as in (2). Using the DIC and MSPE for model selection we obtained the following results:

| | I | | II | | III | |
|---|---|---|---|---|---|---|
| Outcome | MSPE | DIC | MSPE | DIC | MSPE | DIC |

| distribution | | | | | | |
|---|---|---|---|---|---|---|
| Gaussian | 469500 | 1160.7 | 431500 | 1158.8 | 463100 | 1160.3 |

Table 1 Results of fitting a Bayesian linear hierarchical model with Gaussian error with three different random effect assumptions.

Model I has no random effect and is in effect a normal linear model albeit with random coefficients. Model II includes an additive random effect as in equation 1c above. Finally model III is equivalent to a linear mixed model with a measurement error random effect as in equation 2 above.

Both DIC and MSPE suggest that including random effects in the model leads to a better fit. Interestingly, model II appears to be favored over model III. The second model is an overdispersion model, and it must be assumed that for these data the extra variation is reflected in the model fitting rather than error in the DLW covariate. Both DLW and *education level* have significant regression parameters. The level of education has a significant effect on reported energy intake. Women with at least a college degree reported on average 658 kcal/day more than women with less than college degree. Also, women whose energy expenditure was higher reported higher energy intake. Social desirability is not significant under this model. Figure 1 displays a normal quantile plot of the variance-standardized residuals from this model. It is clear that the tails of the distribution are skewed and asymmetry is apparent. However the dominant asymmetry is a heavy left tail and reduced spread in the right tail. This does not correspond closely to a simple log-linear transformation model as proposed by Carroll (Carroll et al.2006). Figure 2 displays the results of fitting models to a simulated data set where a log(y) – log(x) linear relationship is present. The model assumed was $y_i = \beta_0 x_i^{\beta_1} . \exp(e_i)$ which yields a linear model in log(y) and log(x). Figure 2 shows that a log(y)-log(x) transform does fit the data relatively well (left panel). However, a naïve y-x linear model shows a heavy right tail and short left tail i.e. asymmetry. This asymmetry is in the opposite direction from that found in the FFQ_DLW data analyzed here. This certainly challenges the assumption that a log-log transformation can be used within linear mixed models for nutritional measurement error as seems to be commonly assumed.

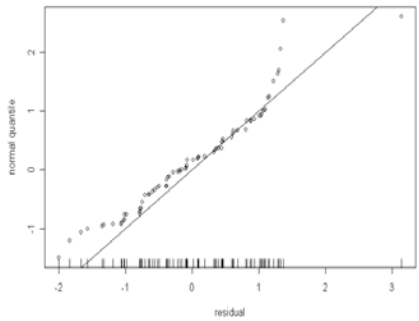

**Figure 1  Standardized predictive residual plot against normal quantiles with associated rug plot:  Normal linear model.**

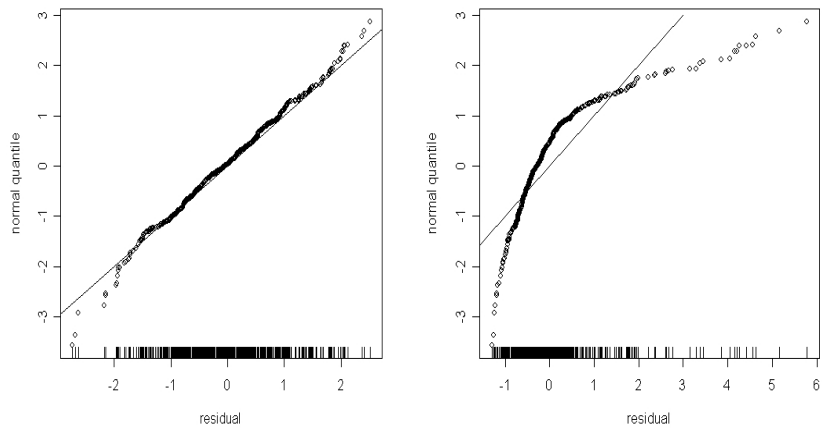

**Figure 2 left panel: normal probability plot of standardized predictive residuals from a log(y)-log(x) linear model when the true model is log(y)-log(x); right panel: normal probability plot of residuals from a naive y-x linear model for the same simulation.**

## 8 Model Extension and Variants

The traditional assumptions of Gaussian linear mixed models have been applied widely in nutritional studies (Carroll et al 1998, Kipnis et al 2003). Their use is supported by the ready availability of software to fit linear mixed models with normal error structures. However these models a priori may be critiqued for a number of reasons:

1) The relation between the dependent variable (a strictly positive quantity) and the covariates (also strictly positive in this case) may be truncated at zero. If the mean response is close to zero a normal distribution would be truncated. This will of course depend on the variance of the response. An asymmetric distribution on the positive real line would be preferred.

2) The overall error in the relationship could be asymmetric and therefore non-normal. This is commonly the case when unobserved heterogeneity could be expected and there is some gross additive error component within the model that is not properly accounted for. In the case of FFQ, the energy calculation is made from weighted sums of macro-nutrients in food servings. These individual self-reported serving errors could be asymmetric and highly correlated. Hence their total effect could easily be conceived to lead to asymmetry in a gross model for this relation, assuming that covariates do not make allowances for these errors. As the covariates are observed at the individual level, it is unlikely that they will compensate for this asymmetry in error.

3) Any measurement error within such models is restricted to zero mean Gaussian forms and this could clearly be atypical of the random measurement error. First, it might be imagined that error could be asymmetric. Self-reporting of food could be biased downwards or upward depending on food type and personality of respondent. This bias could force asymmetry on the responses. Second, errors could be non–linearly related to the variables they are associated with.

4) Transformations of response and covariates could lead to the requirement for strictly positive multiplicative errors both on a linear and exponential scale. Asymmetry will be found when fitting naïve models in un-transformed variables.

5) Repeated measures (Johnson, Herring, and Ibrahim2006) can lead to considerable improvement in the estimation and identification of parameters (random effects ) but may not lead to improved accommodation of the asymmetry in the fundamental relationship

**8.1 Hierarchical Model Alternatives**

To accommodate the skewness in the outcome variable we also considered the following variants of the Gaussian hierarchical model with random effects:

**Log normal variant**

The likelihood and linear predictor are defined as

$$[\log(y_i) \mid \mu_i, \tau_y] \sim N(\mu_i, \tau_y)$$
$$\mu_i = \beta_0 + \beta_1 \log(x_{1i}) + \beta_2 x_{2i} + \beta_3 x_{3i} + \beta_4 (x_{2i} * x_{3i}) + \varepsilon_i \qquad (3)$$

where the random effect $\varepsilon_i$ is shown as an additive effect. A simpler model without random effect has also been examined. These models are labeled as Model II and I in Table 2.

**Prior distribution specification**

The usual specification of prior distributions is made, except that we follow Gelman (2006) and use uniform prior distributions on the standard deviations of regression variance parameters:

$$[\beta_\cdot \mid \tau_{\beta_\cdot}] \sim N(0, \tau_{\beta_\cdot})$$

$$s_\cdot \sim U(0, B_\cdot) \text{ where } s_\cdot = \sqrt{\tau_{\beta_\cdot}} \text{ and } B_* = 1000$$

$$s_y \sim U(0, B_y) \text{ where } s_y = \sqrt{\tau_y} \text{ and } B_e = 1000.$$

We assumed an overdispersed gamma prior distribution for the random effect variance:

$$[\varepsilon_i \mid \tau_\varepsilon] \sim N(0, \tau_\varepsilon)$$

$$\tau_\varepsilon \sim Gamma(\alpha_1, \alpha_2)$$
$$\alpha_1 \sim Gamma(10, 0.1) .$$
$$\alpha_2 \sim Gamma(1, 0.01)$$

Similarly, we also considered a model with random effects inside the regression:

$$[\log(y_i) \mid \mu_i, \tau_y] \sim N(\mu_i, \tau_y)$$
$$\mu_i = \beta_0 + \beta_1 \left( \log(x_{1i}) \varepsilon_i \right) + \beta_2 x_{2i} + \beta_3 x_{3i} + \beta_4 \left( x_{2i} * x_{3i} \right)$$

$$[\beta_\cdot \mid \tau_{\beta_\cdot}] \sim N(0, \tau_{\beta_\cdot})$$

$$s_\cdot \sim U(0, B_\cdot) \text{ where } s_\cdot = \sqrt{\tau_{\beta_\cdot}} \text{ and } B_* = 1000$$

$s_\varphi \sim U(0, B_\varphi)$ where $s_\varphi = \sqrt{\tau_\varphi}$ and $B_c = 1000$

Here we changed the distribution of the random effects to

$[\varepsilon_i \mid \tau_\varepsilon] \sim Gamma(\tau_\varepsilon, \tau_\varepsilon)$ where $\tau_\varepsilon \sim U(0, 500)$.

This model is labeled model III in Table 2.

**Gamma Variant**

Another group of models considered for the outcome variable was based on a Gamma

hierarchical model:

$$[y_i \mid \mu_i, r_y] \sim Gamma(r_y, r), \ r = r_y / \mu_1 \text{ where } E(y_i) = \mu_i$$
$$\mu_i = \beta_0 + \beta_1 (x_{1i} + \varepsilon_i) + \beta_2 x_{2i} + \beta_3 x_{3i} + \beta_4 (x_{2i} * x_{3i}) \qquad (4)$$

Prior distributions are defined as before except that the shape parameter of the gamma likelihood has a gamma prior distribution.

$[\beta. \mid \tau_{\beta.}] \sim N(0, \tau_{\beta.})$

$s. \sim U(0, B.)$ where $s. = \sqrt{\tau_{\beta.}}$ and $B* = 1000$

$[\varepsilon_i \mid \tau_\varepsilon] \sim N(0, \tau_\varepsilon)$

$s_\varepsilon \sim U(0, B_\varepsilon)$ where $s_\varepsilon = \sqrt{\tau_\varepsilon}$ and $B_\varepsilon = 1000$

$r_y \sim Gamma(0.1, 0.001)$

The model defined in (4) is labeled as model III in Table 2. A simpler version without random effect is specified as model I in Table 2. A model with a separate additive random effect is labeled as model II in Table 2.

Table 2 Summary of the results from all fitted models: the significant parameters are listed under the reported MSPE.

|  | I | | II | | III | |
|---|---|---|---|---|---|---|
|  | MSPE | DIC | MSPE | DIC | MSPE | DIC |
| Normal | 469500 $\beta_1$ , $\beta_3$ only | 1160.7 | 431500 $\beta_1$ , $\beta_3$ only | 1158.8 | 463100 $\beta_1$ , $\beta_3$ only | 1160.3 |
| LogNormal | 545200 $\beta_0$ , $\beta_3$ only | 48.45 | 68150 $\beta_0$ only | -190.1 | 507000 $\beta_0$ , $\beta_3$ only | 28.2 |
| Gamma | 490900 $\beta_1$ , $\beta_3$ only | 1159.7 | 463300 $\beta_3$ only | 1157.6 | 490200 $\beta_1$ , $\beta_3$ only | 1160 |

It is clear from Table 2, based on the DIC or MSPE criterion, that the log normal models (with both dependent variable and DLW log-transformed) perform best. However, it is also noticeable that the models with a separate additive random effect term show better performance than those with a term directly associated with the DLW variable. This suggests that there is greater noise in the relation between FFQ and DLW than within DLW itself.

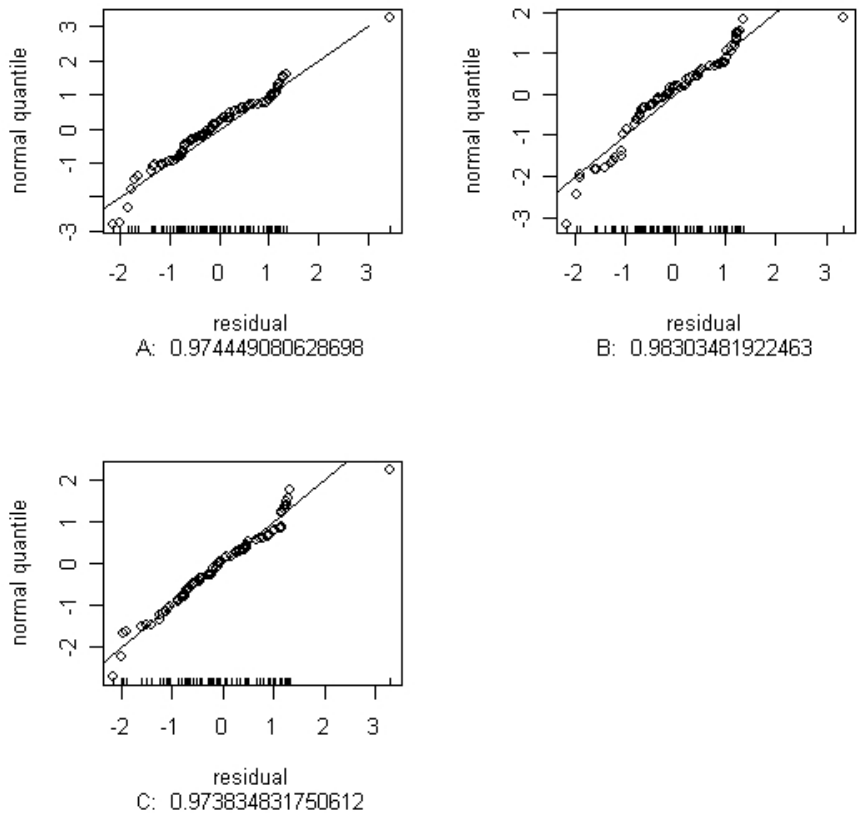

**Figure 3  Normal probability plots: standardized predictive residuals for the models with lognormal. A:**

**model III ( $\beta(\log(x).e)$ );  B: model with  $\beta(\log(x.\varepsilon))$  as a variant to model III; C: model I.**

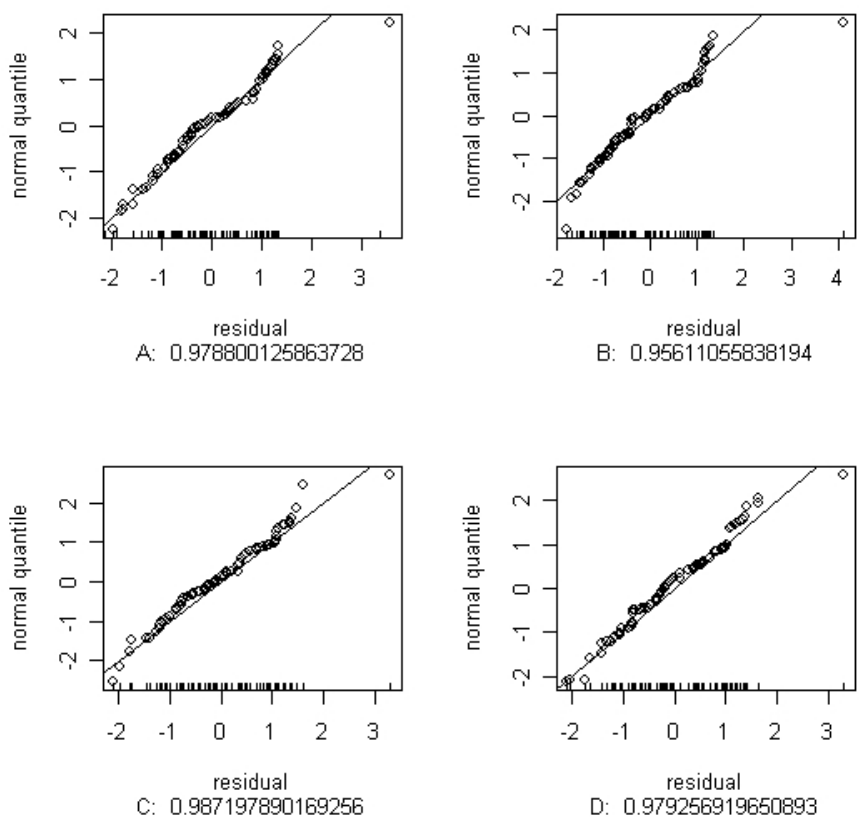

**Figure 4  Normal probability plots of standardized predictive residuals: models with additive random effect and different prior assumptions for normal variance. A:  log(y)-log(x) :  B: log(y)-log(x), C: gamma error; D: gamma error.**

Figure 3 and 4 display the equivalent normal probability plots for the different types of non-Gaussian error models. Beneath each plot is given the empirical Pearson correlation coefficient between the (sorted) residuals and normal quantiles.  In Figure 3, a comparison of log normal models with different random effects is made. In that display it is interesting to note that the models yielding the lowest DIC are those with a multiplicative error (either as ($\beta(\log(x).e)$  or

$\beta (\log (x.\varepsilon)))$ . These models yield the highest correlation although the $\beta (\log (x).e)$ model appears to yield the best linear fit from the appearance of the quantile plot. By contrast the model without random effects has a higher DIC (48.45) but similar quantile plot and correlation to random effect models. Figure 4 displays similar plots but for two different models and prior distributional assumptions for the random effect. As the 'best' models based on DIC and MSPE are models with separate additive random effects, we have examined the effect of variation in prior specification for this type of effect term. The top panel are log (y)-log (x) models with separate random effects (as in model II in Table 3). Plot A has the random effect defined to have a zero mean normal prior distribution but with a uniform hyperprior distribution for the standard deviation of the Gaussian variance (Gelman (2006)). Plot B is the same model but with a gamma prior distribution for the variance as in model (3) above. The bottom panel are gamma error models with separate additive random effects with the same changes in prior distributions. It is interesting to note that while the lowest DIC and MSPE models are the top panel with plot B yielding the 'best' model within Table 2, it is also noteworthy that the highest correlation of residuals to normal quantiles is given by the gamma error models. Indeed across all models in Figure 3 and 4 the lowest DIC model has the lowest correlation. Inspection of the plots suggests that lack-of-fit may be related to a small number of extreme observations. In fact in Figure 4 Plot B two observations appear to be ill-matched to the model fit. Whereas, in Plot D, for example, the model appears to yield residuals that are closer to normal quantiles. While these plots are visual guides to only some aspects of the predictive ability of the models, it is interesting to note that higher DIC or MSPE models can yield residuals that are more regular in distribution. In fact the model with the highest correlation in this regard is a gamma error model with Gelman uniform prior distribution for the random effect variance.

Note that, with the exception of the 'best' log (y)-log (x) model, the gamma error models all have lower MSPE than log normal models. In fact, for the lowest DIC log normal models (the random effect models), the pD parameter which measures the effective number of parameters is negative. This means that the DIC cannot be fully relied on to choose the best model in those cases. Of course the MSPE selects the same overall best log (y)-log (x) model, but the MSPE suggests that the gamma error models perform better than other types of random effect models with log normal errors at the first level of the hierarchy. With respect to the effect of different prior specifications on the random effect terms, this seems to vary depending on the chosen first level error model (log normal or gamma). In the log normal examples, the lowest DIC and MSPE is achieved by the gamma prior model, and it is considerably lower than the Gelman prior specification (DIC: -190.10 compared to -73.58). For the gamma error model the Gelman prior specification yields a marginally better model (DIC: 1157.60 compared to 1159.72).

## 9 Discussion and Conclusions

There are two major conclusions from this work. First, it is clear that the classical formulation of a normal linear mixed model for measurement error does not adequately describe the variety of models we have examined relating FFQ energy to DLW, even with the inclusion of covariates. The departure from this model is particularly marked in the tail areas and demonstrates that there is a strong asymmetry in the error structure. This asymmetry could be due to unobserved confounders that are not accounted for by the random effects included. Essentially, it could be supposed that there could be groupings in the data that when mixed together yield this asymmetry. Overall it can be suggested that transformation (such as log (y)-log (x) for FFQ and DLW) can be used to linearise this relation and produce considerable improvement in model fit.

However it is also clear that such transformation may lead to more serious outlier problems (as demonstrated in Figure 4 and by the reduced correlation coefficients for log (y)-log (x) models) and to an inability to accommodate the variation in relation that is experienced. Hence while such transformations can provide a way to use linear models for FFQ and DLW relations, they lack ability to predict well all the variation found in the relation. Another important issue is the fact that prior specification can have an important impact on the goodness-of-fit of models. Hence, adoption of different prior distributions should be considered within a sensitivity analysis. Given the dependence found here, it is clearly not adequate to make 'standard' assumptions about Gaussian errors and their variances, as is often found in linear mixed model analysis.

Second, for the Energy study, it is clear that models with non-symmetric error structures are appropriate. However, the best fitting models do not include terms in social desirability nor education level, which were thought a priori to make a significant contribution to the relation between FFQ and DLW.

Reference List

1. Berry, S., Carroll, R., and Ruppert, D. 2002. "Bayesian Smoothing and Regression Splines for Measurement Error Problems." *Journal of the American Statistical Association* 97(457):160-169.

2. Brooks, S P. and Gelman, A. 1998. " General Methods for Monitoring Convergence of Iterative Simulations." *Journal of Computational and Graphical Statistics* 7434-55.

3. Carroll, R, Ruppert, D., Stephanski, l., and Crainiceanu, C. 2006. *Measurement Error in Nonlinear Models*. 2nd ed. edited by CRC Press.

4. Carroll, R., Freedman, L., Kipnis, V., and Li, Li. 1998. "A New Class of Measurement-Error Models, With Applications to Dietary Data." *Canadian Journal of Statistics* 26(3):467-77.

5. Gamerman, D. 2002. *Markov Chain Monte Carlo* edited by Chapman&Hall/CRC.

6. Gelfand, A. and Sahu, S. 1999. " Identifiability, Improper Priors, and Gibbs Sampling for Generalized Linear Models." *Journal of the American Statistical Association* 94247-53.

7. Gelman, A. 2006. "Prior Distributions for Variance Parameters in Hierarchical Models." *Bayesian Analysis* 1(2):1-19.

8. Hebert, J., Ebbeling, C. B, Matthews, C. E, Ma, Y., Clemow, L., Hurley, T. G., and Druker, S. 2002. "Systematic Errors in Middle-Aged Women's Estimates of Energy Intake: Comparing Three Self-Report Measures to Total Energy Expenditure From Doubly Labeled Water." *Annals of Epidemiology* 12577-86.

9. Johnson, B. A., Herring, A., and Ibrahim, J. 2006. "Structured Measurement Error in Nutritional Epidemiology: Applications in the Pregnancy, Infection, and Nutrition (PIN) Study." *Journal of the American Statistical Association*(to appear).

10. Kipnis, V., Midthine, D., Freedman, L., Bingham, S., Schatzkin, A., Subar, A., and Carroll, R. 2001. "Empirical Evidence of Correlated Biases in Dietary Assessment Instruments." *American Journal of Epidemiology* 153(4):394-403.

11.  Kipnis, V., Subar, A., Midthune, D., Freedman, L., Ballard-Barbash, R., Troiano, R.,

Bingham, S., Schoeller, D., Schatzkin, A., and Carroll, R. 2003. "Structure of

Dietary Measurement Error: Results of the OPEN Biomarker Study." *American

Journal of Epidemiology* 158(1):14-21.

12.  Rosner, B. and Gore, R. 2001. "Measurement Error Correction in Nutritional Epidemiology

Based on Individual Foods, With Application to the Relation of Diet to Breast

Cancer." *American Journal of Epidemiology* 154(9):827-35.

13.  Spiegelhalter, D., Best, N, Carlin, B, and van der Linde, A. 2002. "Bayesian Measures of

Model Complexity and Fit (With Discussion)." *Journal of the Royal Statistical

Society, Series B* 64(4):583-639.